\documentclass[12pt,reqno]{amsart}
\usepackage{amsmath, amsfonts, amssymb, amsthm}

\usepackage{latexsym}
\def\qed{{\hbadness=10000\hfill\ \vbox{\hrule height.09ex
   \hbox{\vrule width.09ex height1.55ex depth.2ex \kern1.8ex
   \vrule width.09ex height1.55ex depth.2ex}\hrule height.09ex}\break
   \bigskip}}
\newtheorem{theorem} {Theorem} [section]
\newtheorem{proposition} [theorem] {Proposition}
\newtheorem{lemma} [theorem] {Lemma}
\newtheorem{corollary} [theorem] {Corollary}

\theoremstyle{definition}
\newtheorem{definition}[theorem]{Definition}
\newtheorem{notation}[theorem]{Notation}
\newtheorem{remark}[theorem]{Remark}
\def\GL{\operatorname{GL}}
\def\Symp{\operatorname{Sp}}
\def\Polar{\operatorname{W}}
\def\SL{\operatorname{SL}}
\def\PG{\operatorname{PG}}
\def\AG{\operatorname{AG}}

\def\Trace{\operatorname{Trace}}
\def\rad{\operatorname{rad}}

\newcommand{\HH}{{\mathcal{H}}}

\newcommand{\SSS}{{\mathcal{S}}}

\newcommand{\Ff}{{\mathbb F}}

\newcommand{\Ii}{{\mathcal I}}

\newcommand{\llambda}{{\boldsymbol \lambda}}
\newcommand{\Llambda}{{\boldsymbol \Lambda}}

\newcommand{\aalpha}{{\boldsymbol \alpha}}
\newcommand{\bbeta}{{\boldsymbol \beta}}
\newcommand{\sss}{{\mathbf s}}
\newcommand{\rrr}{{\mathbf r}}
\newcommand{\eee}{{\mathbf e}}
\newcommand{\eps}{{\epsilon}}
\newcommand{\abs}[1]{|#1|}

\begin{document}
\title[symplectic groups acting on vectors]
{The permutation action of finite symplectic groups of odd characteristic on their Standard Modules}
\author{David B. Chandler, Peter Sin, Qing Xiang$^*$}

\thanks{$^*$Research supported in part by NSF Grant DMS 0400411.}
\address{Institute of Mathematics, Academia Sinica,
Nangang, Taipei 11529,  Taiwan} \email{chandler@math.udel.edu}
\address{Department of Mathematics, University of Florida, Gainesville, FL 32611,
USA} \email{sin@math.ufl.edu}
\address{Department of Mathematical Sciences, University of Delaware, Newark, DE 19716, USA} \email{xiang@math.udel.edu}

\keywords{Generalized quadrangle, general linear group, $p$-rank,
partial order, symplectic group, symplectic polar space} 

\begin{abstract} Motivated by the incidence problems between points
and flats of a symplectic polar space, we study a large class of
submodules of the space of functions on the
standard module of a finite symplectic group of
odd characteristic. Our structure results on this class of
submodules allow us to determine the $p$-ranks of the incidence
matrices between points and flats of the symplectic polar space. In
particular, we give an explicit formula for the $p$-rank of the
incidence matrix between the points and lines of the symplectic
generalized quadrangle $\Polar(3,q)$, where $q$ is an odd prime
power. Combined with the earlier results of Sastry and Sin on the
2-rank of $\Polar(3,2^t)$, it completes the determination of the
$p$-ranks of $\Polar(3,q)$.
\end{abstract}

\maketitle

\section {Introduction}\label{intro}

Let $k=\Ff_q$ be the finite field of order $q$, where $q=p^t$, $p$
is a prime, and $t$ is a positive integer, and let $V$ be a
$2m$-dimensional vector space over $k$. We denote by $\PG(2m-1 , q)$
the $(2m-1)$-dimensional projective geometry of $V$, and denote by
$P$ the set of points of $\PG(2m-1, q)$. The incidence matrices
between $P$ and flats of $\PG(2m-1, q)$ have been studied
extensively over the past forty years. See for example,
\cite{smith,hamada,Gh,BS,inamdar} for $\Ff_p$-ranks of these matrices, and \cite{sin,csx} for their Smith normal forms. The study of $p$-ranks
of these incidence matrices led the authors of \cite{BS} to
investigate the
submodule lattices of the spaces $k[P]$ and $k[V]$
of $k$-valued functions on $P$ and $V$ respectively,
viewed as permutation modules for the general
linear group $\GL(V)$. The $p$-rank results
can be obtained as a consequence of the description of the submodule
lattice of $k[P]$ (see \cite{BS}). In this paper, we are interested
in certain submatrices of the above mentioned incidence matrices.

We now equip $V$ with a nonsingular alternating bilinear form
$\langle-,-\rangle$. To avoid trivial exceptions, we will assume
that $m\geq 2$ in the rest of this paper. We fix a basis
$e_1,e_2,\ldots ,e_m,f_m,\ldots ,f_1$ and the corresponding
coordinates $x_1,x_2, \ldots ,x_m, y_m,\ldots ,y_1$ so that $\langle
e_i,f_j\rangle=\delta_{ij}$, $\langle e_i,e_j\rangle=0$, and
$\langle f_i,f_j\rangle=0$. The subgroup of $\GL(V)$ leaving
$\langle-,-\rangle$ invariant is the symplectic group $\Symp(V)$.
Let $\Ii_r$ denote the set of totally isotropic $r$-dimensional
subspaces of $V$, where $1\leq r\leq m$. Since $\langle -,-\rangle$
is alternating, we have $\Ii_1=P$, the set of all points of
$\PG(2m-1,q)$. The {\em symplectic polar space} $\Polar(2m-1,q)$ is the
geometry with flats $\Ii_r$, $1\leq r\leq m$. (Here the points of
$\Polar(2m-1,q)$ are the elements of $\Ii_1=P$.) We are interested in
the incidence matrices between points and flats of $\Polar(2m-1,q)$.
More explicitly, for $1\leq r\leq m$, let
\begin{equation}\label{map} \eta_r : k[\Ii_r]\rightarrow k[P]
\end{equation}
be the incidence map sending a totally isotropic $r$-dimensional
subspace of $V$ to its characteristic function in $P$. We are
interested in the images of the maps $\eta_r$. These incidence
problems concerning $\Polar(2m-1,q)$ lead naturally to the
study of $k\Symp(V)$-submodules of $k[P]$ and $k[V]$.

Our main results, under the assumption that $q$ is odd,
are as follows. We
will define a special basis of $k[V]$ (see
Definition~\ref{sympbasisdef}), whose elements are called {\it symplectic
basis functions}. Our main theorems describe the submodule structure
of the $k\Symp(V)$-module generated by an arbitrary symplectic basis
function. The $k\Symp(V)$-module $k[P]$ can be viewed as a direct
summand of $k[V]$, and the class of submodules of $k[V]$ described
above includes the images of $\eta_r$, $1\leq r\leq m$. We then
obtain $p$-rank formulas for the incidence matrices between
$\Ii_1=P$ and $\Ii_r$ from the structure results on the submodules
mentioned above. In the case where $m=2$, we obtain a particularly
nice $p$-rank formula, which we will describe below in some detail.

For convenience, let $A_{1,r}^m(q)$ be a $(0,1)$-matrix with rows
indexed by the elements $Y$ of $\Ii_r$ and columns indexed by the
elements $Z$ of $P$, and with the $(Y,Z)$ entry equal to 1 if and
only if $Z\subseteq Y$. We consider the case where $m=2$ (and $r=2$)
in particular. In this case, the symplectic polar space $\Polar(3,q)$
is a classical generalized quadrangle (GQ) \cite{Ti,paynethas}, whose
points are all the points of $\PG(3,q)$, and whose lines are the
totally isotropic 2-dimensional subspaces of $V$. When $q=2^t$,
Sastry and Sin \cite{ss} gave the following formula for the 2-rank
of $A_{1,2}^2(q)$.
\begin{equation}\label{2rank}
{\rm rank}_2(A_{1,2}^2(2^t))=1+\left(\frac {1+\sqrt
{17}}{2}\right)^{2t}+\left(\frac {1-\sqrt {17}}{2}\right)^{2t}.
\end{equation}
In the case where $q=p$ is an odd prime, de Caen and Moorhouse
\cite{moorhouse} determined the $p$-rank of $A_{1,2}^2(p)$, which was
later generalized by the second author \cite{sinsymp}, giving the
$p$-ranks of $A_{1,r}^m(p)$, where $1\leq r\leq m$, $p$ is an odd
prime, and $m$ is not necessarily 2. In this paper, we obtain the
following formula for the $p$-rank of $A_{1,2}^2(p^t)$, $p$ an odd
prime, as a corollary of our submodule structure results.

\begin{theorem}\label{prank}
 Let $p$ be an odd prime and let $t\geq 1$ be an integer. Then the $p$-rank of
 $A_{1,2}^2(p^t)$ is equal to
 $$1+\alpha_1^t+\alpha_2^t,$$
where
\begin{equation}\label{oddrk}
\alpha_1,\alpha_2 =
\frac{p(p+1)^2}{4}\pm\frac{p(p+1)(p-1)}{12}\sqrt{17}.
\end{equation}
\end{theorem}

We remark that in (\ref{oddrk}), if we simply set $p=2$, then we
actually obtain (\ref{2rank}), but the two results require different
proofs.

The paper is organized as follows. In Section~\ref{Action}, we will review the
results in \cite{BS} concerning the $\GL(V)$-submodule lattice of
$k[V]$.
The submodule lattice has a combinatorial description in terms of
certain partially ordered sets $\HH$ and $\HH[d]$.
(See subsection~\ref{HHtypes}
 below.) For the moment, we will just consider $\HH$, which is
associated with the nontrivial summand $Y_P$ of $k[P]$. The module
$Y_P$ has a special basis, and to each basis element there is
an associated element of $\HH$ called its $\HH$-type, giving a
surjective map from the basis to $\HH$. It was proved in \cite{BS}
that for each $\sss\in\HH$, the set of basis elements whose
$\HH$-types are $\leq \sss$ span a $k\GL(V)$-submodule $Y(\sss)$ of
$Y_P$ with the property that $Y(\sss)$ has a unique maximal
submodule. Furthermore, every submodule of $Y_P$ is a sum of
submodules of the form $Y(\sss)$.

On the representation-theoretic side, the main goal of this paper is
to construct analogues of these objects adapted to the action of
$\Symp(V)$. In order to do so, it is necessary first to look
deeper into the $k\GL(V)$-structure of $k[V]$. By considering its
multiplicative structure as a $k\GL(V)$-algebra, we derive tensor
product factorizations of certain subquotients of $k[V]$ which will
be needed in our later constructions. These new results concerning
$\GL(V)$ are also included in Section~\ref{Action}. In Section~\ref{posets}, we define
posets $\SSS$ and $\SSS[d]$ whose elements are pairs $(\sss,\eps)$,
with $\sss$ in $\HH$ (or $\HH[d]$) and $\eps$ a certain
``signature''. In Section~\ref{Action4}, we define a special basis of $k[V]$.
Just as in the $\GL(V)$ case, a certain subset of this basis spans
$Y_P$ and there is surjection from this subset to $\SSS$. For
$(\sss,\eps)\in \SSS$, let $Y(\sss,\eps)$ be the $k$-subspace
spanned by the basis elements of $Y_P$ which map into the ideal in
$\SSS$ determined by $(\sss,\eps)$. In Section~\ref{Submodules} we prove that
$Y(\sss,\eps)$ is a $k\Symp(V)$-submodule of $Y_P$, and our main
technical result, that $Y(\sss,\eps)$ has a unique maximal
submodule. Unlike the $k\GL(V)$-submodules,
not every $\Symp(V)$-submodule of $Y_P$ is the sum of submodules of
the form $Y(\sss,\eps)$.  The reason is a fundamental difference between the
two cases.  As a $k\GL(V)$-module, $Y_P$ is multiplicity-free---that
is, no two composition factors are isomorphic---while the
 $k\Symp(V)$-module is not. Nevertheless, the portion of the entire
$k\Symp(V)$-submodule lattice generated by the submodules
$Y(\sss,\eps)$ is sufficiently rich for our applications. In
Section~\ref{dimensions}, we apply the results of Section~\ref{Submodules}
 to ${\rm Im}(\eta_r)$, the
images of the incidence maps $\eta_r$ defined in (\ref{map}).
In this way, we obtain a summation formula for the
$p$-rank of the incidence matrix $A_{1,m}^m(p^t)$, where $p$ is odd.
In particular, we give a proof of Theorem~\ref{prank}.

\section{Action of $\GL(V)$ on $k[V]$}\label{Action}

Throughout Sections~\ref{Action} through \ref{Submodules} of the paper, we assume that $p$
is an odd prime, $k=\Ff_q$, $V$ is a $2m$-dimensional vector space
over $k$, and $q=p^t$, $t>1$. The assumption that $t>1$ is mainly
for notational convenience, and is only seriously used in
Lemma~\ref{shift} and Lemma~\ref{digitprojector}. We shall need to
apply some of the results of \cite{BS}.

The results in \cite[Theorems A, B, C]{BS}
give a simple and complete description
of the  $k\GL(V)$-submodule structure of the space $k[V]$
of $k$-valued functions on a finite vector space $V$.
Let $k[X_1, X_2,\ldots,X_{2m}]$ denote the polynomial ring, in $2m$
variables. Since every function on $V$ is given by a polynomial in
the $2m$ coordinates $x_i$, the map $X_i\mapsto x_i$ defines a
surjective $k$-algebra homomorphism $k[X_1,
X_2,\ldots,X_{2m}]\rightarrow k[V]$, with kernel generated by the
elements $X_i^q-X_i$. Furthermore, this map is simply the coordinate
description  of the following  canonical map. The polynomial ring is
isomorphic  to the symmetric algebra $S(V^*)$ of the dual space of
$V$, so we have a natural evaluation map $S(V^*)\rightarrow k[V]$.
This canonical description makes it clear that the map is
equivariant with respect to the natural actions of $\GL(V)$ on these
spaces. A basis for $k[V]$ is obtained by taking monomials in $2m$
coordinates $x_i$ such that the degree in each variable is at most
$q-1$. We will call these the {\it basis monomials} of $k[V]$.

The space $k[V]$ has the structure of a $\mathbb Z/(q-1)\mathbb
Z-$graded $\GL(V)$-algebra, where the grading is given by the
characters of the center, the scalar multiplications, isomorphic to
$k^\times$. Thus,
$$
k[V]=\oplus_{[d]\in \mathbb Z/(q-1)\mathbb Z}A[d],
$$
where $\mu\in k^\times$ acts on the component $A[d]$
as $\mu^{[d]}$. The component $A[d]$ has basis consisting of the basis
monomials in which  the total degree is in the residue class $[d]$.

\subsection{Types and $\HH$-types}\label{HHtypes}
We now recall the definitions of two $t$-tuples associated
with each basis monomial. Let
\begin{equation}\label{monomial}
f=\prod_{i=1}^{2m}x_i^{b_i}=\prod_{j=0}^{t-1}\prod_{i=1}^{2m}(x_i^{a_{ij}})^{p^j},
\end{equation}
be a basis monomial,
where $b_i=\sum_{j=0}^{t-1}a_{ij}p^j$ and $0\leq a_{ij}\leq p-1$. Let $\lambda_j=\sum_{i=1}^{2m}a_{ij}$. The $t$-tuple
$\llambda=(\lambda_0,\ldots,\lambda_{t-1})$ is called the {\it type} of $f$.
The set of all types of monomials is denoted by $\Llambda$.

Let $d$ be the integer between $0$ and $q-2$ which is congruent to
the total degree $\sum_ib_i=\sum_j\lambda_jp^j$ modulo $q-1$, and
let $(d_0,\dots, d_{t-1})$  be the  $t$-tuple of $p$-adic digits of
$d$.

In \cite{BS}, there is another $t$-tuple
associated with each basis monomial, which
we will call its \emph{$\HH$-type.}
If $[d]\neq[0]$ this tuple will lie in the set
\begin{equation*}
\HH[d]=\lbrace\sss=(s_0,\ldots,s_{t-1})\mid \forall j,
0\leq s_j\leq {2m-1}, 0\leq ps_{j+1}-s_j +d_j\leq 2m(p-1)\rbrace,
\end{equation*}
and if $d=0$, it will belong to the set $\HH[0]=\HH\cup\{(0,0,\ldots
0), (2m,2m,\ldots 2m)\}$, where
\begin{equation*}
\HH=\lbrace \sss=(s_0,s_1,\dots,s_{t-1})\mid \forall j, 1\leq s_j\leq 2m-1,\
0\leq ps_{j+1}-s_j \leq 2m(p-1) \rbrace.
\end{equation*}
The $\HH$-type $\sss$ of $f$ is uniquely determined by the type via
the equations
\begin{equation*}
\lambda_j=ps_{j+1}-s_j+d_j, \quad 0\leq j\leq t-1.
\end{equation*}
Moreover, these equations determine a bijection between the set
$\Llambda$ of types of basis monomials and the union of the sets
$\HH[d]$, $0\leq d\leq q-2$. We will consider the sets $\HH[d]$ and
$\HH$ as partially ordered sets under their natural order induced by
the product order on $t$-tuples of natural numbers.

\begin{notation}\label{omit}
We will be considering many objects indexed by $\HH$-types.
To indicate that the corresponding $\HH$-type belongs
to $\HH[d]$, a decoration $[d]$ will be used. In the case
$[d]=0$, we will most often be interested in the case where
the $\HH$-type is in $\HH$. In this case, we adopt the convention
of omitting $[0]$ from the notation.
\end{notation}

\subsection{Composition factors}
The types, or equivalently the $\HH$-types parametrize the
composition factors of $k[V]$ in the following sense.
Except for the existence of two trivial direct summands in $A[0]$,
the $k\GL(V)$-module $k[V]$ is multiplicity-free. We can
associate to each $\HH$-type $\sss\in\HH[d]$ a composition factor,
which we shall denote by
$L(\sss)[d]$, such that these simple modules are all
nonisomorphic except that $L((0,\dots,0))[0]\cong L((2m,\ldots,2m))[0]
\cong k$. The simple modules $L(\sss)[d]$ occur as subquotients
of $k[V]$ in the following way.
For $\rrr\in\HH[d]$
let $Y(\rrr)[d]$ be the span of all basis
monomials with $\HH$-types in $\HH[d]_\rrr=\{\rrr'\in\HH[d]\mid \rrr'\leq\rrr\}$. By \cite{BS}, if $[d]\neq[0]$ then $Y(\rrr)[d]$ is a
$k\GL(V)$-submodule of $A[d]$ with a unique maximal submodule
and such that the quotient by the
maximal  submodule is isomorphic to $L(\rrr)[d]$.
In the case $[d]=[0]$, for each $\sss\in\HH$ we let $Y(\sss)$,
be the subspace spanned by monomials of $\HH$-types
in $\HH_\sss=\{\sss'\in\HH\mid \sss'\leq\sss\}$, and
similarly, $Y(\sss)$ has a unique simple quotient, isomorphic
to $L(\sss):=L(\sss)[0]$ (by the notational convention above).

The isomorphism type of the simple module $L(\sss)[d]$ is most
easily  described in terms of the corresponding type
$(\lambda_0,\ldots,\lambda_{t-1})\in \Llambda$. Let $S^\lambda$ be
the degree $\lambda$ component in the truncated polynomial ring
$k[X_1, X_2,\ldots,X_{2m}]/(X_i^p; 1\leq i\leq 2m)$. Here $\lambda$
ranges from $0$ to $2m(p-1)$. Note that the dimension of $S^\lambda$
is
\begin{equation}\label{dimen}
d_\lambda=\sum_{j=0}^{\lfloor \lambda/p\rfloor}(-1)^j {2m\choose
j}{2m-1+\lambda-jp \choose 2m-1}.
\end{equation}
The simple module $L(\sss)[d]$ is isomorphic to the twisted tensor
product
\begin{equation}\label{tensorproduct}
S^{\lambda_0}\otimes(S^{\lambda_1})^{(p)}\otimes\cdots\otimes
(S^{\lambda_{t-1}})^{(p^{t-1})}.
\end{equation}

\begin{remark}\label{truncatedbasis}
Note that each module $(S^\lambda)^{(p^j)}$ is itself isomorphic to
a composition factor $L(\sss)[p^j\lambda]$ of  $k[V]$, corresponding
to the type $\llambda$ with $\lambda_j=\lambda$ and all other
components zero. Let us be more precise about this identification.
From the definition, we may view $(S^\lambda)^{(p^j)}$ as the degree
$\lambda$ component of the truncated polynomial ring in the
variables $X_i^{p^j}$.  In $k[X_1,\ldots,X_{2m}]$ we consider the
set of $p^j$-th powers of monomials of total degree $\lambda$ and
with the degree of each variable between $0$ and $p-1$. This set
maps injectively into the truncated polynomial ring in the variables
$X_i^{p^j}$ and the images form a basis for $(S^\lambda)^{(p^j)}$.
The images of the same monomials in $k[V]$ are basis monomials of
type $\llambda$. Hence they lie in $Y(\sss)[p^j\lambda]$ and they map
bijectively to a basis of the simple quotient $L(\sss)[p^j\lambda]$.
Later on, when we  abuse notation slightly and speak of
$(S^\lambda)^{(p^j)}$ as having a basis consisting of images of
basis monomials of type $\llambda$, the exact meaning will always be
as we have just described.
\end{remark}

\subsection{Submodule structure}
The reason for considering $\HH$-types is that they allow a simple
description of the submodule structure of the $k\GL(V)$-modules
$A[d]$.  Suppose first that $[d]=[0]$. The space $A[0]$
has a trivial  direct summand spanned by the characteristic function
of $\{0\}$, which is the kernel of the natural map
$A[0]\rightarrow k[P]$.
The basis monomials with types in $\HH[0]$, excluding the
type $(2m(p-1),2m(p-1),\ldots,2m(p-1))$ span a complementary
direct summand, which maps isomorphically onto $k[P]$.
We have
\begin{equation}\label{splittingofk[P]}
A[0]\cong k\oplus k[P]=k\oplus k\oplus Y_P,
\end{equation}
where $Y_P$ is the kernel of the map $k[P]\rightarrow
k$, $f\mapsto |P|^{-1}\sum_{Q\in P}f(Q)$. The $k\GL(V)$ module $Y_P$
is an indecomposable module whose composition factors are
parametrized by $\HH$. The \cite[Theorem A]{BS} states that given
any $k\GL(V)$-submodule of $Y_P$, the set of its composition factors
is an ideal in the partially ordered set $\HH$ and that this
correspondence is an order isomorphism from the submodule lattice of
$Y_P$  to the lattice of ideals in $\HH$. For a submodule $A\leq
Y_P$, let $\HH(A)\subseteq\HH$ denote the ideal of $\HH$-types of
its composition factors.

For $[d]\neq[0]$, the set $\HH[d]$ parametrizes the composition
factors of $A[d]$ and we have a similar order isomorphism
from the submodule lattice of $A[d]$ to the lattice
of ideals in $\HH[d]$, with its natural partial order
\cite[Theorem C]{BS}.
Let $\HH[d](A)$ denote the ideal of the submodule $A\leq A[d]$.

Assume now that $M$ is a subquotient of $A[d]$ with no trivial
submodules. This condition is just a convenient way of saying
that in the case $[d]=0$ we assume $M$ is a subquotient of $Y_P$
(so that its set of $\HH$-types is well defined).
Then there are submodules $B\leq C$ of $A[d]$
with no trivial submodules such that $M=C/B$. Thus, if
$[d]\neq 0$, the composition
factors of $M$ correspond to the set
$\HH[d](C)\setminus\HH[d](B)$, which is
a difference of  ideals in $\HH[d]$, while if $[d]=[0]$
the composition factors of $M$ correspond to
$\HH(C)\setminus\HH(B)$, a difference of ideals in $\HH$.

The submodules of $A[d]$ and $Y_P$ can also be described in terms of
basis monomials \cite[Theorem B]{BS}. Any submodule of $A[d]$
($[d]\neq[0]$) or of $Y_P$ has a basis consisting of the basis
monomials which it contains.  Moreover, the $\HH$-types of
these basis monomials are precisely the $\HH$-types of the
composition factors of the submodule. Furthermore, in any
composition series, the images of the monomials of a fixed
$\HH$-type form a basis of the composition factor of that
$\HH$-type. (These statements are not quite true of $A[0]$,
because of the two trivial summands.)

%
\subsection{$\GL(V)$-algebra structure}
Multiplication in $k[V]$ is pointwise multiplication of functions
and it is $\GL(V)$-equivariant, giving  $k\GL(V)$-homomorphisms
\begin{equation*}
A[d]\otimes A[d']\rightarrow A[d+d'],\quad \text{for}\quad[d], [d']\in\mathbb Z/(q-1)\mathbb Z.
\end{equation*}

\begin{lemma} Let $\llambda=(\lambda_0,\ldots,
\lambda_{t-1})\in \Llambda$ correspond to the $\HH$-type
$\rrr=(r_0,\ldots,r_{t-1})\in \HH[d]$. Let
$[d^*]=[d-\lambda_{t-1}p^{t-1}]$ and let the $\HH$-type
$\rrr^*=(r^*_0,\ldots,r^*_{t-1})\in\HH[d^*]$ correspond to the type
$\llambda^*=(\lambda_0,\ldots,\lambda_{t-2},0)$. Then

\begin{equation}
\rrr^*=\rrr+\eee,
\end{equation}
where the $t$-tuple $\eee$ of integers depends only on $[d]$ and $\lambda_{t-1}$.
\end{lemma}
\begin{proof} The lemma follows directly from the definitions of
$\rrr$ and $\rrr^*$. Let $d_j$ and $d^*_j$ be the $p$-adic digits of
the least nonnegative residues in $[d]$ and $[d^*]$ respectively.
Then by definition,
$$
\lambda_j=pr_{j+1}-r_j+d_j,\qquad \lambda^*_j=pr^*_{j+1}-r^*_j+d^*_j;
$$
so for $0\leq i\leq t-1$,
$$
(q-1)r_i=\sum_{j=0}^{t-1}(\lambda_j-d_j)p^{(j-i)},
$$
where the exponent $(j-i)$ is taken to be the least nonnegative
residue modulo $t$. The lemma follows by comparing this formula with
the similar one for $r^*_i$, remembering that
$\lambda^*_i=\lambda_i$ for $0\leq i\leq t-2$ and that $[d^*]$ is
determined by $[d]$ and $\lambda_{t-1}$.
\end{proof}
\begin{corollary} Let $\mathcal T\subseteq \Llambda$ be
a set of types whose $(t-1)$-th entries are all equal to
$\lambda_{t-1}$.
Let $\mathcal T^*$ be the set of types
obtained from $\mathcal T$ by replacing $\lambda_{t-1}$
by zero in the $(t-1)$-th entry. Let $\mathcal X$
and $\mathcal X^*$ be the corresponding subsets of
$\HH$-types in $\HH[d]$ and $\HH[d^*]$ respectively, with
the induced orderings. The following hold:
\begin{enumerate}
\item[(i)] The bijection $\mathcal T\rightarrow
\mathcal T^*$ sending $\llambda$ to $\llambda^*$
induces an order isomorphism from $\mathcal X$ to $\mathcal X^*$.
\item[(ii)] $\mathcal X$ is a difference of ideals
of $\HH[d]$ if and only if  $\mathcal X^*$
is a difference of ideals of $\HH[d^*]$.
\end{enumerate}
\end{corollary}
\begin{proof} Both follow from the previous lemma; for (ii)
we note that a subset of a finite partially ordered set is a
difference of ideals if and only if it satisfies the ``intermediate
value'' condition that for any two elements in the subset, all
elements in between them are also in the subset.
\end{proof}

\begin{theorem} Let $M$ be a $k\GL(V)$-subquotient of $A[d]$ with
no trivial submodules and let $\mathcal X$ denote the set of $\HH$-types of its composition factors in $\HH[d]$.
Suppose that for some $j\in \mathbb Z/t\mathbb Z$,
all tuples in $\mathcal X$ have the same $r_j$ and also the same
entries $r_{j+1}$. Let $\lambda_j=pr_{j+1}-r_j+d_j$.
Let $\mathcal T\subseteq \Llambda$ be the set of types
corresponding to $\mathcal X$ and $\mathcal T^*$ be the set of types
obtained from $\mathcal T$ by replacing $\lambda_j$
by zero in the $j$-th entry.

Then in the $k\GL(V)$-submodule
$P$ of $A[d-p^j\lambda_j]$ generated by all monomials
with types in $\mathcal T^*$, the $k$-subspace $Q$ spanned by
monomials whose types are not in $\mathcal T^*$ is a
$k\GL(V)$-submodule. Let $N=P/Q$. Then
\begin{equation}
M\cong N\otimes (S^{\lambda_j})^{(p^j)}.
\end{equation}
Moreover, the types of $N$ are obtained by replacing
$\lambda_j$ by $0$ in the types of $M$.
\end{theorem}

\begin{proof}
Note that in the case $[d]=[0]$ our hypothesis implies that  $M$ is a
subquotient of $Y_P$, so the set of $\HH$-types of its composition
factors is well-defined. By Galois conjugation, we may assume
$j=t-1$. By \cite[Theorem A]{BS} $\mathcal X$ is a difference of
ideals of $\HH[d]$.
  Let $\mathcal T$ be the set of
types $\llambda=(\lambda_0,\ldots\lambda_{t-1})$ corresponding to
$\mathcal X$. By hypothesis, the entry
$\lambda_{t-1}=d_{t-1}+pr_0-r_{t-1}$ is the same for every type in
$\mathcal T$. Then, by the previous corollary, the set $\mathcal
X^*\subseteq\HH[d^*]$, whose types form the set $\mathcal T^*$ of
types obtained from $\mathcal T$ by replacing $\lambda_{t-1}$ by $0$
in the $(t-1)$-th entry, is a difference of ideals in $\HH[d^*]$,
where $[d^*]=[d-\lambda_{t-1}p^{t-1}]$. Let $P\leq A[d^*]$ be the
$k\GL(V)$-submodule generated by all monomials of types in $\mathcal
T^*$.  Then by \cite{BS} there exists a $k\GL(V)$-submodule $Q\leq P$
such that $Q$ has as basis all the monomials of $P$ whose types are
not in $\mathcal T^*$, and $N=P/Q$ is a $k\GL(V)$-module with basis
consisting of the bijective images of all monomials of type $\mathcal T^*$.
Likewise, $M$ has a basis
consisting of images of all monomials whose types lie in $\mathcal
T$. In exactly the same way, the $p^{t-1}$-th powers of all
monomials of degree $\lambda_{t-1}$ form a basis of a
$k\GL(V)$-subquotient $S$ of $A[\lambda_{t-1}p^{t-1}]$ with
$S\cong(S^{\lambda_{t-1}})^{(p^{t-1})}$ as $k\GL(V)$-modules.

It is clear that if we multiply each monomial of type $\mathcal T^*$
by the $p^{t-1}$-th power of each monomial of degree $\lambda_{t-1}$,
we obtain each monomial of type $\mathcal T$ exactly once.
Therefore the multiplication map
$A[d^*]\otimes A[\lambda_{t-1}p^{t-1}]\rightarrow A[d]$
induces a bijection of the subquotients
\begin{equation}
\label{tensoriso}
N\otimes S \cong M.
\end{equation}
Since the multiplication map is a map of $k\GL(V)$
modules, the map (\ref{tensoriso}) is a $k\GL(V)$-isomorphism.
\end{proof}

\begin{remark}\label{interpretiso} Let us interpret this tensor factorization
in terms of a function $f\in A[d]$ which maps to a nonzero element
$\overline f$ of $M$. Assume that $f$ can be written as a product
$f=f'{f_j}^{p^j}$, where the monomials of $f'\in A[d-p^j\lambda_j]$
have types in $\mathcal T^*$ and those of ${f_j}^{p^j}\in
A[p^j\lambda_j]$ are of type $(0,\ldots,0,\lambda_j,0,\ldots,0)$.
Then under the isomorphism of the theorem, $\overline f$ is mapped
to $\overline f'\otimes\overline {{f_j}^{p^j}}$, where $\overline
f'$ is the image of $f'$ in the subquotient $N$ of
$A[d-p^j\lambda_j]$ and $\overline{{f_j}^{p^j}}$ is the image of
${f_j}^{p^j}$ in the simple subquotient $S$ of
$A[\lambda_{j}p^{j}]$.
\end{remark}

\subsection{The modules $Y(\sss)[d]_j$ and $Y(\sss)_j$}
We will consider certain quotients of $Y(\sss)[d]$ and $Y(\sss)$.
Let $\mathcal X\subset\HH[d]_\sss$ be the subset of tuples having
$j$-th and $(j+1)$-th entries equal to $s_j$ and $s_{j+1}$
respectively and $\lambda_j=m(p-1)$. It is clear  that $\mathcal X$
is the difference of the ideal $\HH[d]_\sss$ and an ideal of
$\HH[d]$, since it satisfies the ``intermediate value'' condition;
so $\mathcal X$ is the set of tuples of a $k\GL(V)$-quotient
$\overline Y(\sss)[d]_j$ of $Y(\sss)[d]$. Moreover, in the case
$[d]=[0]$, we have $\mathcal X\subseteq\HH$ and so $\overline
Y(\sss)[0]_j$ is actually  a quotient of $Y(\sss)$. The following is
immediate from the theorem above.

\begin{lemma}\label{factorize} There is a $k\GL(V)$-module
$B_j$ such that
\begin{equation*}
\overline Y(\sss)[d]_j\cong  B_j\otimes (S^{m(p-1)})^{(p^j)}.
\end{equation*}
\end{lemma}

\vspace{0.1in}

\section{The posets $\SSS$ and $\SSS[d]$}\label{posets}
\begin{definition}
For $\llambda\in\Llambda$, let $\sss$ be the corresponding $\HH$-type
in $\HH[d]$. Set
$$J(\sss)=\{j \mid 0\leq j\leq t-1,\;\lambda_j=m(p-1)\}.$$
For any $\sss,\sss'\in\HH[d]$, let
$Z(\sss,\sss')=\{j \mid s'_j=s_j,\ s'_{j+1}=s_{j+1},\
\lambda_j=m(p-1)\}$. We define
$$\SSS[d]=\{(\sss,\eps)\mid \sss\in \HH[d],\ \eps\subseteq J(\sss)\}.$$
In the case $[d]=[0]$, we also define
$$\SSS=\{(\sss,\eps)\mid \sss\in \HH, \eps\subseteq J(\sss)\}.$$

We define $(\sss',\eps')\leq (\sss,\eps)$ if and only if $\sss'\leq\sss$
and $\eps\cap Z(\sss',\sss)=\eps'\cap Z(\sss',\sss)$. It is not difficult
to check that this defines a partial order on $\SSS[d]$
and $\SSS$; for transitivity
one notes that if $\sss''\leq\sss'\leq\sss$ then
$Z(\sss'',\sss)=Z(\sss'',\sss')\cap Z(\sss',\sss)$.
\end{definition}

Since each $\sss\in\HH$ or $\rrr\in\HH[d]$ corresponds to a type
$\lambda\in \Llambda$, we can also talk about signed types
$(\llambda,\eps)$ corresponding to elements of $\SSS$ or $\SSS[d]$.

\section{Action of $\Symp(V)$ on $k[V]$}\label{Action4}
We now equip $V$ with a nonsingular alternating bilinear form
$\langle-,-\rangle$, with the basis $e_1,e_2,\ldots ,e_m,f_m,\ldots
,f_1$ and the corresponding coordinates $x_1,x_2, \ldots ,x_m,
y_m,\ldots ,y_1$ as given in Section~\ref{intro}. Accordingly, we view
$S(V^*)$ as the polynomial ring  generated by
``symplectic indeterminates'', $X_1$,\ldots $X_m$, $Y_m$,\ldots, $Y_1$.

We will consider the
submodule structures of $k[V]$, $A[d]$, and $Y_P$, under the action of
$\Symp(V)$. First let us recall the known facts about composition
factors (cf. \cite{SZ,La}). We would like to know how a
$\GL(V)$ composition factor (\ref{tensorproduct}) decomposes upon
restriction to $\Symp(V)$. The modules $S^\lambda$,
$0\leq\lambda\leq 2m(p-1)$, all remain simple except when
$\lambda=m(p-1)$, in which case we have
\begin{equation}\label{splitS}
S^{m(p-1)}=S^+\oplus S^-.
\end{equation}
Here, $S^+$ and $S^-$ are simple $k\Symp(V)$-modules, and
\begin{equation}\label{dim+-}
{\rm dim}(S^+)=(d_{(p-1)m}+p^m)/2,\quad {\rm
dim}(S^-)=(d_{(p-1)m}-p^m)/2.
\end{equation}

We can describe $S^+$ and $S^-$ as follows.

To avoid cumbersome notation involving $X_1$,\ldots, $X_m$,
$Y_m$,\ldots, $Y_1$, we will use multi-index notation $X^{\aalpha}
Y^{\bbeta}$ for monomials, where ${\aalpha}=(a_1,\ldots,a_m)$ and
${\bbeta}=(b_1,\ldots,b_m)$, $0\leq a_i,b_i\leq p-1$. Further, for
any multi-index ${\bbeta}$, we define
$\abs{{\bbeta}}=\sum_{i=1}^mb_i,\ {\bbeta}!=\prod_{i=1}^mb_i!$,
and $\overline {\bbeta}=(p-1-b_1,\ldots,p-1-b_m)$, and similarly define
$\abs{\aalpha}$ and $\aalpha!$. We will denote
the images of monomials in the simple module $S^{m(p-1)}$ using
bars. Then \cite{La} the map
\begin{equation}
\tau:S^{m(p-1)}\rightarrow S^{m(p-1)}, \quad \overline X^{\aalpha}
\overline Y^{\bbeta}\mapsto (-1)^{\abs{{\bbeta}}}
{\aalpha}!{\bbeta}!\overline X^{\overline{\bbeta}}\overline
Y^{\overline{\aalpha}}
\end{equation}
is a $k\Symp(V)$-homomorphism with $\tau^2=1$.

The modules $S^+$ and $S^-$ are the eigenspaces of $\tau$ for the
eigenvalues $(-1)^m$ and $(-1)^{m+1}$ respectively. By
Remark~\ref{truncatedbasis} the space $S^{m(p-1)}$ can be viewed as
having a basis of images of basis monomials of $k[V]$. From this
point of view, the eigenspaces $S^+$ and $S^-$ have bases consisting
images of basis monomials of $k[V]$ of the form
\begin{equation}\label{basis1}
x^{\aalpha} y^{\overline{\aalpha}}
\end{equation}
and of sums and differences
\begin{equation}\label{basis2}
x^{\aalpha} y^{\bbeta}\pm
(-1)^{\abs{{\bbeta}}+m}{\aalpha}!{\bbeta}!x^{\overline{\bbeta}}y^{\overline{\aalpha}}
\end{equation}
of monomials, for ${\aalpha}\neq\overline{\bbeta}$. The images of
the monomials (\ref{basis1}) together with those of the form
(\ref{basis2}) with a ``$+$'' sign form a basis of $S^+$ and those
with a ``$-$'' sign form a basis of $S^-$.

\begin{definition}\label{sympbasisdef}
We will now define a new basis of $k[V]$, whose elements we will
call {\it symplectic basis functions}. We will first define
the symplectic basis functions of type $\llambda$. Then we will
take the union of these sets of functions over all $\llambda$.
The symplectic basis functions of type $\llambda$ will
be certain functions of the form
\begin{equation}\label{bf} f=f_0f_1^p\cdots f_{t-1}^{p^{t-1}}. \end{equation}
where each $f_j$, which we will call the $j$-th digit of $f$, is
either a basis monomial or binomial of $k[V]$ of degree $\lambda_j$.
We will now describe the allowable forms of the $j$-th digit; then
the set of functions $f$, all of whose digits are allowable, will be
the set of symplectic basis functions of type $\llambda$. If
$\lambda_j\ne (p-1)m$, then $f_j$ can be any basis monomial of
degree $\lambda_j$ in which the degree in each variable
is at most $p-1$. If $\lambda_j=(p-1)m$, then $f_j$ can be any
function of the form (\ref{basis1}) or (\ref{basis2}).
\end{definition}

Clearly by restricting the types for the symplectic basis functions
we can obtain bases for $A[d]$,  and $Y_P$.

\begin{definition}
To each symplectic basis function of $k[V]$ we associate a pair
$(\sss,\eps)\in\SSS[d]$ for some $[d]\in \mathbb Z/(q-1)\mathbb Z$,
 as follows. If $f$ is of type
$\llambda$, then $\sss$ is the corresponding $\HH$-type. The set
$\eps\subseteq J(\sss)$, called the \emph{signature,} is defined to be
the set of $j\in J(\sss)$ for which the image of the $j$-th digit
$f_j$ of $f$ in $S^{m(p-1)}$ belongs to $S^+$.
\end{definition}

From (\ref{tensorproduct}) and (\ref{splitS}), it is clear that the
$k\Symp(V)$-composition factors of $k[V]$ are given by their types,
together with the additional choice of signs for each $j$ with
$\lambda_j=m(p-1)$. In terms of $\HH$-types, we see that each
$\HH$-type gives a $k\GL(V)$-composition factor and then the choice
of signs determines the simple $k\Symp(V)$ composition factor of
this simple $k\GL(V)$-module. In this way, the elements of $\SSS$
label the $k\Symp(V)$-composition factors of $Y_P$, and those of
$\SSS[d]$, $[d]\neq [0]$ label the $k\Symp(V)$-composition factors
of $A[d]$. However it should be noted that different elements of
$\SSS$ or $\SSS[d]$ can label isomorphic composition factors, due to
the fact that $S^\lambda\cong S^{2m(p-1)-\lambda}$ as
$k\Symp(V)$-modules. We will use $L(\sss,\eps)[d]$ to denote the
simple $k\Symp(V)$-submodule of $L(\sss)[d]$ where we take the $+$
summand for each $j\in\eps$ and the $-$ summand for each $j\in
J(\sss)\setminus\eps$. When $\sss\in \HH$, we may use the simpler
notation $L(\sss,\eps)$.

It follows from the definitions that the set of symplectic basis
functions of $\HH$-type $\sss\in \HH[d]$ and signature $\eps$ maps
bijectively under the natural map $Y(\sss)[d]\rightarrow L(\sss)[d]$
to a basis of $L(\sss,\eps)[d]$. We will also call $\eps$ the
signature of $L(\sss,\eps)[d]$.

The following statement is an immediate consequence of Lemma~\ref{factorize} and the
decomposition of $S^{m(p-1)}$ just discussed.
\begin{theorem}\label{plusminus}
As $k\Symp(V)$-modules, we have
\begin{equation*}
\overline Y(\sss)[d]_j\cong (B_j\otimes (S^+)^{(p^j)})\oplus(B_j\otimes (S^-)^{(p^j)}).
\end{equation*}
\end{theorem}

\vspace{0.1in}

\section{The submodules $Y(\sss,\eps)[d]$}\label{Submodules}

Let $Y(\sss)[d]_j^+$ be the preimage in $Y(\sss)[d]$
 of the $+$ component of
$\overline Y(\sss)[d]_j$ in Theorem~\ref{plusminus} and let
$Y(\sss)[d]_j^-$ be the preimage in $Y(\sss)[d]$ of the $-$
component. For $\eps\subseteq J(\sss)$, let
\begin{equation}\label{Y(s,e)}
Y(\sss,\eps)[d]=\bigcap_{j\in\eps}Y(\sss)[d]^+_j\quad{\textstyle\bigcap}\quad\bigcap_{j\in J(\sss)\setminus \eps}Y(\sss)[d]^-_j.
\end{equation}
Thus, $Y(\sss,\epsilon)[d]$  is a  $k\Symp(V)$-submodule of
$Y(\sss)[d]$.

\begin{lemma}\label{signedbasis} Let $(\sss,\eps)\in\SSS[d]$.
Then $Y(\sss,\eps)[d]$ has a basis consisting of all the symplectic
basis functions with signed $\HH$-types
$(\sss',\eps')\leq(\sss,\eps)$.
\end{lemma}
\begin{proof}
Suppose $(\sss',\eps')\in\SSS[d]$ with
$(\sss',\eps')\leq(\sss,\eps)$. and let $f$ be a symplectic basis
function of signed type $(\sss',\eps')$. We will show first that
$f\in Y(\sss,\eps)[d]$. Write $f=f_0f_1^p\cdots f_{t-1}^{p^{t-1}}$
as the product of its digits raised to the appropriate powers. Let
$j\in J(\sss)$.  We must show that $f\in Y(\sss)[d]_j^+$ if
$j\in\eps$ and $f\in Y(\sss)[d]_j^-$ if $j\in J(\sss)\setminus\eps$.

If $f$ maps to zero in $\overline Y(\sss)[d]_j$ then it is clear
from the definitions that $f\in Y(\sss,\eps)$[d]. So we may assume
that $f$ has nonzero image $\overline f\in\overline Y(\sss)[d]_j$.
According to Remark~\ref{interpretiso}, under the isomorphism of
Lemma~\ref{factorize}, $\overline f$ is  mapped to $\overline
f'\otimes\overline f_j^{p^j}$, where $\overline f_j^{p^j}$ is the
image of $f_j^{p^j}$ in $(S^{m(p-1)})^{(p^j)}$  and $\overline f'$
is the image in $B_j$ of the product of the other factors of $f$.
Thus, since $\overline f\neq 0$, we must have $j\in Z(\sss,\sss')$. 
From the definition of $\tau$ and the assumption that $f_j$ has an
allowable form, we see that $\overline f_j^{p^j}$ is an eigenvector
of the endomorphism of $(S^{m(p-1)})^{(p^j)}$ induced by $\tau$.
Therefore $\overline f$ is an eigenvector of the endomorphism of
$\overline Y(\sss)[d]_j$ induced by $\tau$ via the tensor
factorization of Lemma~\ref{factorize}, and will belong to either
the $+$ or $-$ part of the decomposition given in
Theorem~\ref{plusminus}. More precisely, $\overline f$ will be in
the $+$ part if $j\in \eps'$ and in the $-$ part if $j\in
J(\sss')\setminus\eps'$. But we already have $j\in Z(\sss,\sss')$
and since $(\sss',\eps')\leq(\sss,\eps)$, we have $\eps\cap
Z(\sss,\sss')=\eps'\cap Z(\sss,\sss')$. Thus $\overline f$ is in the
$+$ part if $j\in\eps$ and in the $-$ part if $j\notin\eps$. We have
proved $f\in Y(\sss,\eps)[d]$.

Now we must prove that
$Y(\sss,\eps)[d]$ is spanned by the symplectic basis functions with
signed types $(\sss',\eps')\leq (\sss,\eps)$. Since we know that
$Y(\sss)[d]$ has a basis consisting of all symplectic basis
functions with $\HH$-types $\sss'\leq\sss$, it suffices to prove
that no linear combination
\begin{equation}\label{gs}
\sum_i c_i g_i,
\end{equation}
with nonzero scalars $c_i$, of symplectic basis functions whose
signed types $(\sss_i,\eps_i)$ satisfy $\sss_i\leq\sss$ but
$(\sss_i,\eps_i)\nleq(\sss,\eps)$, can belong to $Y(\sss,\eps)[d]$.
Consider the function $g_1$. There must exist $j\in Z(\sss,\sss_1)$
which belongs to $\eps$ but not $\eps_1$, or \emph{vice versa}.  We will assume $j\in\eps$, as the case $j\in J(\sss)\setminus\eps$ is similar.
 We can rewrite (\ref{gs})
as
\begin{equation*}
\sum_{i\in I}c_i g_i +\sum_{r\notin I}c_r g_r,
\end{equation*}
where $I$ is the set of indices $i$ for which $j\in Z(\sss,\sss_i)$.
Under the map $Y(\sss)[d]\rightarrow\overline Y(\sss)[d]_j$, the set
$\{g_i\mid i\in I\}$ is mapped to a linearly independent set, while
the elements $g_r$ with $r\notin I$ are mapped to zero. The reason is
 that $\overline Y(\sss)[d]_j$ corresponds to the set of
$\HH$-types $\sss'\leq\sss$ for which $s'_j=s_j$, $s'_{j+1}=s_{j+1}$,
and $\lambda'_j=m(p-1)$. Therefore, the image in $\overline
Y(\sss)[d]_j$ of $\sum_{i\in I}c_ig_i$ is a sum of linearly
independent eigenvectors for the endomorphism induced by $\tau$. At
least one of the terms, namely the image $g_1$, has the opposite
eigenvalue to that prescribed by $\epsilon$. The conclusion is that the
image of $\sum_{i\in I}c_ig_i$ in $\overline Y(\sss)[d]_j$ cannot be
in the $+$ component of $\overline Y(\sss)[d]_j$ as given in
Theorem~\ref{plusminus}. Therefore $\sum_ic_ig_i$ cannot belong to
$Y(\sss,\eps)[d]$. The proof is complete.
\end{proof}

It is obvious from Lemma~\ref{signedbasis} that
$Y(\sss',\eps')[d]\leq Y(\sss,\eps)[d]$ if and only if
$(\sss',\eps')\leq(\sss,\eps)$. We define $Y_<(\sss,\eps)[d]$ to be
the kernel of the natural map of $Y(\sss,\eps)[d]\rightarrow
L(\sss,\eps)[d]$, or equivalently, the sum of all
$Y(\sss',\eps')[d]$ with $(\sss',\eps')\lneq (\sss,\eps)$.

\begin{remark}\label{digitsum}
For $\sss\in\HH[d]$, we define its {\it digit sum} by
$\abs{\sss}=\sum_{j=0}^{t-1}s_j$. It is not hard to see that if
$(\sss',\eps')\lneq (\sss,\eps)$ then there exists $(\sss'',\eps'')$
such that $\abs{\sss''}= \abs{\sss}-1$ and $(\sss',\eps')\leq
(\sss'',\eps'')\leq (\sss,\eps)$; so we also have
\begin{equation}\label{Y<}
Y_<(\sss,\eps)[d]=\sum_{\begin{smallmatrix}
(\sss',\eps')\leq(\sss,\eps)\\
\abs{\sss'}=\abs{\sss}-1\end{smallmatrix}} Y(\sss',\eps')[d].
\end{equation}
\end{remark}

\subsection{Submodule structure of $Y(\sss,\eps)[d]$}
\begin{lemma} \label{multfree}
Let $\sss\in\HH[d]$. Then no two $k\Symp(V)$ composition factors of
$\oplus_{\sss'} L(\sss')[d]$ are isomorphic, where the sum runs over
all $\sss'$ which are immediately below $\sss$.
\end{lemma}
\begin{proof} Let $\sss'$ and $\sss''\in\HH[d]$ be immediately
below $\sss$. It is clear that no two simple $k\Symp(V)$-submodules
of $L(\sss')[d]$ are isomorphic; so we must consider the case where
some $L(\sss',\eps')[d]$ is isomorphic to some
$L(\sss'',\eps'')[d]$. Now these simple modules have the form of
twisted tensor products, which, by Steinberg's Tensor Product
Theorem can be isomorphic only if the corresponding tensor factors
are isomorphic. Thus, the above isomorphism can only happen if
$L(\sss')[d]$ and $L(\sss'')[d]$ are isomorphic as
$k\Symp(V)$-modules, which means that for each $j$ we must have
either $\lambda'_j=\lambda''_j$ or $2m(p-1)-\lambda''_j$. Let
$\sss=(s_0,\ldots,s_{t-1})$, with similar notation for $\sss'$ and
$\sss''$. By Galois conjugation we may assume without loss that
$s'_0=s_0-1$ and $s''_k=s_k-1$ for some $k\neq 0$. Suppose first
$t>2$ and $k\neq 1$. Then $\lambda'_0=\lambda_0+1$ and
$\lambda''_0=\lambda_0$,  so that the above condition cannot hold. If
$k=1$, then by considering $\lambda'_1=\lambda_1$ and
$\lambda''_1=\lambda_1+1$, we reach the same conclusion. Finally we
must consider the case $t=2$. Then the above condition forces
$2\lambda_0=2\lambda_1=(2m+1)(p-1)$. Therefore $s_0=s_1$ and so
$\lambda_0=(p-1)s_0$. Dividing the previous equation by $(p-1)$
yields the desired contradiction.
\end{proof}
\vspace{-0.2in} Fix $[d]\in\mathbb Z/(q-1)\mathbb Z$
and$(\sss,\eps)\in\SSS[d]$.

Let $\mathcal Z$ be the set of elements of $\SSS[d]$, which are
immediately below $(\sss,\eps)$. Let
$R=\sum_{(\sss'',\eps'')\in\mathcal Z}Y_<(\sss'',\eps'')[d]$. Then
\begin{equation*}
Y_<(\sss,\eps)[d]/R\cong\oplus_{(\sss'',\eps'')\in\mathcal Z}
(Y(\sss'',\eps'')[d]/Y_<(\sss'',\eps'')[d])
\end{equation*}
is a multiplicity-free semisimple module by Lemma~\ref{multfree}.
Fix $(\sss',\eps')\in\mathcal Z$ and let
\begin{equation*}
K(\sss',\eps')=Y_<(\sss',\eps')[d]+\sum_
{\begin{smallmatrix}(\sss''\eps'')\in\mathcal Z\\
                    (\sss'',\eps'')\neq(\sss',\eps')\end{smallmatrix}}Y(\sss'',\eps'')[d],
\end{equation*}
and let $U=Y(\sss,\eps)[d]/Y_<(\sss,\eps)[d]$.  Then we have a short exact sequence
\begin{equation}\label{ses1}
0\rightarrow
(Y(\sss',\eps')[d]+K(\sss',\eps'))/K(\sss',\eps')\rightarrow
Y(\sss,\eps)[d]/K(\sss',\eps')\rightarrow U
\rightarrow 0.
\end{equation}
which is an extension of $L(\sss,\eps)[d]$ by $L(\sss',\eps')[d]$.

We will show that the short exact sequence (\ref{ses1}) does not
split. To do so, we need to introduce some shift operators
(elements in the group ring $k\Symp(V)$). The $p$-adic version of
these shift operators was used extensively in \cite{csx}. Here we
are using the finite field version of these operators. For $\mu\in
k^{\times}$, we use $g_{\mu}$ to denote the symplectic transvection
sending $x_1$ to $x_1+\mu y_1$ and fixing all other coordinates.

\begin{lemma}\label{shift}
For $0\leq j\leq t-1$ and $1\leq \ell\leq p-1$, let
\begin{equation}\label{shiftoper}
g_{\ell}(j) = \sum_{\mu\in k^{\times}} \mu^{\ell p^j}g_{\mu^{-1}}
\in k\Symp(V).
\end{equation}
Given any basis monomial $f=x_1^{a_1}y_1^{b_1}\cdots
x_m^{a_m}y_m^{b_m}$ of $k[V]$, we have
$$g_{\ell}(j) f= \left\{\begin{array}{ll}
    0, & \textrm{\emph{if the $j$-th digit of $a_1$ is less than $\ell$;}} \\
    -{a_1 \choose \ell} x_1^{a_1-\ell p^j}y_1^{b_1+ \ell p^j}x_2^{a_2}y_2^{b_2}\cdots x_m^{a_m}y_m^{b_m}, & \mathrm{otherwise.}
            \end{array}\right. $$
\end{lemma}

\begin{proof} We first prove the lemma for $j=0$. If $a_1=0$, then
clearly we have $g_\ell (0)f=0$. So we assume that $a_1 > 0$.
\begin{eqnarray*}
g_{\ell}(0) f &= & \sum_{\mu\in k^{\times}} \mu^{\ell} (x_1 + \mu^{-1} y_1)^{a_1}y_1^{b_1}x_2^{a_2}y_2^{b_2}\cdots x_m^{a_m}y_m^{b_m}\\
&= &\left(\sum_{\mu\in k^{\times}} \mu^{\ell} (x_1^{a_1} + {a_1 \choose 1}\mu^{-1}x_1^{a_1 -1}y_1+{a_1\choose 2}\mu^{-2} x_1^{a_1 -2}y_1^2 +\cdots)\right) y_1^{b_1}x_2^{a_2}y_2^{b_2}\cdots x_m^{a_m}y_m^{b_m}\\
&= & -{a_1\choose \ell} x_1^{a_1-\ell}y_1^{b_1
+\ell}x_2^{a_2}y_2^{b_2}\cdots x_m^{a_m}y_m^{b_m}
\end{eqnarray*}
By a classical theorem of Lucas \cite{lucas}, ${a_1 \choose
\ell}\equiv 0$ (mod $p$) if the $0^\mathrm{th}$ digit of $a_1$ (in
the base $p$ expansion of $a_1$) is less than $\ell$, proving the
lemma for $j=0$.

The general case follows from the $j=0$ case by using the Frobenius
automorphism.
\end{proof}

We let $h_{\ell}(j)$ denote the group ring element analogous to
$g_\ell(j)$, but with the roles of $x_1$ and $y_1$ exchanged, so
that this element shifts $\ell p^j$ from the exponent of $y_1$ to
that of $x_1$.

\begin{lemma}\label{digitprojector}
For each pair of integers $(\alpha,\beta),\ 0\le \alpha,\beta\le
p-1$, and each $j,\ 0\le j\le t-1$, there is a group ring element
$g_{\alpha,\beta}(j)\in k\Symp(V)$ such that for any basis monomial
\begin{equation}\label{basef}
f=\prod_{i=1}^m {x_i}^{a_i}{y_i}^{b_i}
\end{equation}
 of $k[V]$, where
$a_i=\sum_{k=0}^{t-1}a_{ik}p^k\ \mathrm{and}\  b_i=\sum_{k=0}^{t-1}b_{ik}p^k$,\
$0\leq a_{ik},b_{ik}\leq p-1$,
$$g_{\alpha,\beta}(j)f= \left\{\begin{array}{ll}
    f, & \mathrm{if}\ a_{1j}=\alpha\ \mathrm{and}\ b_{1j}=\beta,\
\mathrm{or}\ a_{1j}=p-1-\beta\ \mathrm{and}\ b_{1j}=p-1-\alpha; \\
    0, & \mathrm{otherwise.}
            \end{array}\right. $$
\end{lemma}

\begin{proof} It  suffices to prove the lemma in the case where
$\alpha +\beta\leq p-1$. The reason is that for any pair
$(\alpha,\beta),\ 0\le \alpha,\beta\le p-1$, with $\alpha +\beta>
p-1$, the ``complementary'' pair $(p-1-\beta, p-1-\alpha)$ has sum of
entries equal to $2(p-1)-(\alpha +\beta)$, which is $<p-1$, and
$g_{p-1-\beta, p-1-\alpha}(j)$ will be the required element. We will
only give the proof for the $j=0$ case. The other cases are the
same. We use induction on $\alpha +\beta$.

First assume that $\alpha +\beta=p-1$. Using Lemma~\ref{shift}, we
define
\begin{equation}\label{operator}
g_{\alpha,\beta}(0)=-{p-1\choose\beta}^{-1}g_\beta(0)h_{p-1}(0)
g_\alpha(0).\end{equation}
We claim that $g_{\alpha,\beta}(0)$ has
the required action on the basis monomials of $k[V]$. It can be
seen as follows. Let $f$
be a
basis monomial of $k[V]$ as in (\ref{basef}).
We
first assume that $a_{10}+b_{10}\leq p-1$. By Lemma~\ref{shift},
$$g_{\alpha}(0)f=\left\{\begin{array}{ll}
    -{a_{10}\choose \alpha}x_1^{a_1-\alpha}y_1^{b_1+\alpha}x_2^{a_2}y_2^{b_2}\cdots x_m^{a_m}y_m^{b_m}, & \mbox{if $a_{10}\geq \alpha$;} \\
    0, & \mathrm{otherwise.}
    \end{array}\right. $$
Next,
$$h_{p-1}(0)(g_{\alpha}(0)f)=\left\{\begin{array}{ll}
    {a_{10}\choose \alpha}x_1^{a_1-\alpha +p-1}y_1^{0}x_2^{a_2}y_2^{b_2}\cdots x_m^{a_m}y_m^{b_m}, & \mbox{if $a_{10}\geq \alpha$ and $b_{10}+\alpha =p-1$;} \\
    0, & \mathrm{otherwise.}
    \end{array}\right. $$
Note that $\alpha +\beta =p-1$, and we have $b_{10}+\alpha =p-1$ if and
only if $b_{10}=\beta$. So
$$h_{p-1}(0)(g_{\alpha}(0)f)=\left\{\begin{array}{ll}
    {a_{10}\choose \alpha}x_1^{a_1-b_{10}}y_1^{0}x_2^{a_2}y_2^{b_2}\cdots x_m^{a_m}y_m^{b_m}, & \mbox{if $a_{10}\geq \alpha$ and $b_{10}=\beta$;} \\
    0, & \mathrm{otherwise.}
    \end{array}\right. $$
Finally,
$$g_{\beta}(0)(h_{p-1}(0)g_{\alpha}(0)f)=\left\{\begin{array}{ll}
    -{p-1 \choose \beta}{a_{10}\choose \alpha}x_1^{a_1}y_1^{b_1}x_2^{a_2}y_2^{b_2}\cdots x_m^{a_m}y_m^{b_m}, & \mbox{if $a_{10}\geq \alpha$ and $b_{10}=\beta$;} \\
    0, & \mathrm{otherwise.}
    \end{array}\right. $$
Note that under the assumptions $\alpha +\beta =p-1$ and
$a_{10}+b_{10}\leq p-1$, the condition that $a_{10}\geq \alpha$ and
$b_{10}=\beta$ means exactly $a_{10}=\alpha$ and $b_{10}=\beta$.
We have shown that the group ring element defined in (\ref{operator})
has the desired action on those $f$ with $a_{10}+b_{10}\leq p-1$.

For the purpose of this lemma, two monomials $f$ and $f'$ can be considered
``complementary," if  $f=\prod_{i=1}^mx_i^{a_i}y_i^{b_i}$ and
$f'=\prod_{i=1}^mx_i^{a'_i}y_i^{b'_i}$ with
$a_{10}+b'_{10}=a'_{10}+b_{10}=p-1$.  The elements we are constructing should act the same on $f'$ and on $f$.  In particular, the element
(\ref{operator}) acts the same on $f'$ as it does on $f$. That is,
$$g_{\alpha,\beta}(j)f'= \left\{\begin{array}{ll}
    f', & \mbox{if $a_{10}=\alpha$ and $b_{10}=\beta$;} \\
    0, & \mathrm{otherwise.}
            \end{array}\right. $$
The
analysis is quite similar to the above.  (One needs to take
extra care when there is a carry from the $0$-th digit to the first
digit, such as in the case where $p-1-a_{10}+\alpha \geq p$. We omit the details.) With
this observation, we see that the group ring element defined in
(\ref{operator}) also has the desired action on those $f$ with
$a_{10}+b_{10}> p-1$, proving the base case where
$\alpha +\beta =p-1$.

For a general pair $(\alpha,\beta)$ with $\alpha +\beta <p-1$, by
induction hypothesis, we may assume that for all those pairs
$(\gamma,\delta)$, $0\leq \gamma, \delta\leq p-1$, with $\alpha
+\beta< \gamma +\delta <p$, we have found $g_{\gamma, \delta}(0)$
with the desired property. We define
\begin{equation*}\label{genoperator}
g_{\alpha,\beta}(0)=-{\alpha+\beta\choose\beta}^{-1}g_\beta(0)h_{(\alpha+\beta)}(0)
g_\alpha(0)
\prod_{\alpha+\beta<\gamma+\delta<p}\left(1-g_{\gamma,\delta}(0)\right).
\end{equation*}
Again we claim that this $g_{\alpha,\beta}(0)$ has the required
action on the basis monomials as given in (\ref{basef}).
Clearly, if $\alpha
+\beta<a_{10}+b_{10}<2(p-1)-(\alpha +\beta)$, then $f$ will be
annihilated by
$\prod_{\alpha+\beta<\gamma+\delta<p}\left(1-g_{\gamma,\delta}(0)\right)$.
So we only need to consider the action of $g_{\alpha,\beta}(0)$ on
those $f$ with
\begin{equation*}\label{survive}
 a_{10}+b_{10}\leq \alpha +\beta<p-1
\quad \mathrm{or}\quad p-1<2(p-1)-(\alpha +\beta)\leq a_{10}+b_{10}.
\end{equation*}
It is clear that
$\prod_{\alpha+\beta<\gamma+\delta<p}\left(1-g_{\gamma,\delta}(0)\right)$
acts on such basis monomials as the identity, and we only need to
consider the action of
$-{\alpha+\beta\choose\beta}^{-1}g_\beta(0)h_{(\alpha+\beta)}(0)
g_\alpha(0)$ on these monomials.

Now if $a_{10}+b_{10}<p$, an
analysis similar to that in the $\alpha +\beta =p-1$ case shows that
$$-{\alpha+\beta\choose\beta}^{-1}g_\beta(0)h_{(\alpha+\beta)}(0)
g_\alpha(0)(f)=\left\{\begin{array}{ll}
    f, & \mbox{if $a_{10}=\alpha$ and $b_{10}=\beta$;} \\
    0, & \mathrm{otherwise,}
    \end{array}\right. $$
and in the complementary case $a_{10}+b_{10}\ge p-1,$
$$-{\alpha+\beta\choose\beta}^{-1}g_\beta(0)h_{(\alpha+\beta)}(0)
g_\alpha(0)(f')=\left\{\begin{array}{ll}
    f', & \mbox{if $a'_{10}+\beta=p-1$ and $\alpha+b'_{10}=p-1$;} \\
    0, & \mathrm{otherwise.}
    \end{array}\right. $$
The proof is complete.
\end{proof}

\begin{lemma}\label{nonsplit}
Assume that $(\sss,\eps)$ is not a minimal element of $\SSS[d]$.
If $[d]=[0]$ we assume in addition that $(\sss,\eps)\in\SSS$
and is not minimal in $\SSS$. Then the short exact sequence
(\ref{ses1}) does not split.
\end{lemma}
\begin{proof}
We will choose a particular element $f\in Y(\sss,\eps)[d]$
with nonzero image in
$$U=Y(\sss,\eps)[d]/Y_<(\sss,\eps)[d]$$
and show that
if $f^*\in Y(\sss,\eps)[d]$ is any element
with the same image in $U$, 
then as a $k\Symp(V)$-module, $Y(\sss,\eps)[d]/K(\sss',\eps')$ is
generated by the image
of $f^*$ in $Y(\sss,\eps)[d]/K(\sss',\eps')$.
Let us first fix some notation. Since $\sss'$ is immediately
below $\sss$, there is a unique index $j+1$
where these tuples differ and $s'_{j+1}=s_{j+1}-1$.
We let $\llambda$ and $\llambda'$ be the corresponding types.

The element $f$ is chosen to be a certain symplectic basis function
of signed $\HH$-type $(\sss,\eps)$. Since $(s_0,\ldots
,s_j,s_{j+1}-1,s_{j+2},\ldots ,s_{t-1})\in\HH$, we have
$\lambda_j\geq p$ and $\lambda_{j+1}<2m(p-1)$. Therefore we may
choose $f$ such that the $j$-th digit of the exponent of $x_1$ is
least 1 and the $j$-th digit of the exponent of $y_1$ is equal to
$p-1$. We can also require that the $(j+1)$-th digit of the exponent
of $y_1$ be less than $p-1$, and further, if
$\lambda_{j+1}=m(p-1)-1$, that the $(j+1)$-th digits of the exponents
of $x_1$ and $y_1$ be 0. Let $f_j$ denote the $j$-th digit of $f$.

Let $e=f^*-f\in Y(\sss',\eps')[d]+K(\sss',\eps')$. From the
definition of symplectic basis functions $f^*$ has the form
\begin{equation*}
f^*=f_0f_1^p\cdots (x_1^{a}y_1^{p-1}\cdots )^{p^j}f_{j+1}^{p^{j+1}}\cdots f_{t-1}^{p^{t-1}} + e
\end{equation*}
or
\begin{equation*}
f^*=f_0f_1^p\cdots (x_1^{a}y_1^{p-1}\cdots\pm cx_1^0y_1^{p-1-a}\cdots )^{p^j}f_{j+1}^{p^{j+1}}\cdots f_{t-1}^{p^{t-1}} + e,
\end{equation*}
where $a\geq 1$, $c$ represents the product of factorials as in
(\ref{basis2}), and $f_{j+1}$ could be a monomial or another term of
the same form as in (\ref{basis1}) or (\ref{basis2}).

Now we apply the group ring element $g_{a,p-1}(j)$ from
Lemma~\ref{digitprojector} to $f^*$. It annihilates all but those
monomials appearing in $e$ with the same $j$-th digits of the exponents of $x_1$ and $y_1$ as those of $f$ or the complementary $j$-th digits,
$0$ and $p-1-a$, respectively. Next we apply the shift operator
$g_1$ from Lemma~\ref{shift} which shifts $p^j$ from the exponent of
$x_1$ to that of $y_1$. The results are
\begin{equation}\label{goodelt}
\begin{aligned}
&g_1(f_0f_1^p\cdots (x_1^{a}y_1^{p-1}\cdots )^{p^j}f_{j+1}^{p^{j+1}}\cdots f_{t-1}^{p^{t-1}})\\
&=f_0f_1^p\cdots (x_1^{a-1}y_1^{0}\cdots )^{p^j}(f_{j+1}y_1)^{p^{j+1}}\cdots f_{t-1}^{p^{t-1}},\end{aligned}
\end{equation}
which is of type $\sss'$, and
\begin{equation}
g_1(f_0f_1^p\cdots (x_1^{0}y_1^{p-1-a}\cdots )^{p^j}f_{j+1}^{p^{j+1}}\cdots f_{t-1}^{p^{t-1}})=0.
\end{equation}
Note that if $f'$ is any other monomial in $f^*$ belonging to
$(Y(\sss',\eps')[d]+K(\sss',\eps'))$, not annihilated by
$g_{a,p-1}(j)$, then $g_1(f')\in K(\sss',\eps')$ because the
$\HH$-type of $g_1(f')$ is obtained by subtracting 1 from the
$(j+1)$-th entry of the $\HH$-type of $f'$. Now we have produced an
element (\ref{goodelt}) of $(Y(\sss',\eps')[d]+K(\sss',\eps'))$ and
we must show that it is not zero, modulo $K(\sss',\eps')$, or in
other words, that when this element is expressed in symplectic basis
functions, necessarily of type $\sss'$, some symplectic basis
function of signature $\eps'$ appears with nonzero coefficient.
Since $(\sss',\eps')\leq(\sss,\eps)$ and the tuples $\sss$ and
$\sss'$ differ only in the $(j+1)$-st digit, we have
$Z(\sss,\sss')\subseteq J(\sss')\subseteq
Z(\sss,\sss')\cup\{j,j+1\}$. We consider several possibilities. If
$J(\sss')=Z(\sss,\sss')$ then all basis functions involved in
(\ref{goodelt}) are of signature $\eps'$. If $j\in J(\sss')$, that
is to say $\lambda'_j=m(p-1)$, then in (\ref{goodelt}) the $j$-th
digit of the exponent of $x_1$ is at most $p-2$, and that of $y_1$
is 0. The monomial is not of the form (\ref{basis1}) and therefore
can be written as the sum of a (nonzero) $S^+$ term and a (nonzero)
$S^-$ term. Similarly, if $\lambda_{j+1}'=m(p-1)$, we have taken
care that the $(j+1)$-th digits of the exponents of $x_1$ and $y_1$
are $0$ and $1$, and the monomial is not of the form (\ref{basis1}).
Thus, for all four possibilities for $J(\sss')$ and for any
$\eps'\subseteq J(\sss')$ such that $\eps'\cap
Z(\sss,\sss')=\eps\cap Z(\sss,\sss')$, the element in
(\ref{goodelt}) involves a symplectic basis function with signature
$\eps'$. The proof is now complete.
\end{proof}

\begin{remark}
If $[d]=[0]$, the assumption of Lemma~\ref{nonsplit} is equivalent
to $\sss\neq (0,\ldots,0)$, $(2m,\ldots,2m)$ or $(1,\ldots,1)$.
\end{remark}

We recall that the {\it radical} of a module is the intersection
of its maximal submodules. Let $\rad M$ denote the radical of
a $k\Symp(V)$-module $M$. It is the largest submodule of $M$
such that the quotient is semisimple.

\begin{theorem}\label{submodule}
\emph{(i)} If $[d]\neq[0]$ then $Y_<(\sss,\eps)[d]$ is the unique maximal
$k\Symp(V)$-submodule of $Y(\sss,\eps)[d]$.

\emph{(ii)} For $(\sss,\eps)\in\SSS$,
$Y_<(\sss,\eps)$ is the unique maximal
$k\Symp(V)$-submodule of $Y(\sss,\eps)$.
\end{theorem}

\begin{proof}
We may assume in both parts that $(\sss,\eps)$ satisfies the
hypotheses of Lemma~\ref{nonsplit}, for otherwise $Y(\sss,\eps)[d]$
in (i) and $Y(\sss,\eps)$ in (ii) are simple modules and there is
nothing to prove. We will only give the argument for (ii), since the
proof of (i) is formally identical. The assertion of the theorem can
be restated as $\rad Y(\sss,\eps)=Y_<(\sss,\eps)$. We proceed by
induction on the partial order of $\SSS$, the result being clear for
the minimal element. Let $f\in Y(\sss,\eps)\setminus Y_<(\sss,\eps)$
and let $Y_f$ be the $k\Symp(V)$ submodule generated by $f$. By
Lemma~\ref{nonsplit} the sequence (\ref{ses1}) does not split.
Therefore $Y_f$ contains an element of
$Y(\sss^*,\eps^*)+K(\sss^*,\eps^*)$ which has nonzero image in
$$
(Y(\sss^*,\eps^*)+K(\sss^*,\eps^*))/K(\sss^*,\eps^*)\cong
Y(\sss^*,\eps^*)/Y_<(\sss^*,\eps^*).
$$
Thus, $(Y_f+R)/R$ has $L(\sss^*,\eps^*)$ as a composition factor.
Since $(\sss^*,\eps^*)$ was an arbitrary element of $\mathcal Z$ and
since $Y_<(\sss,\eps)/R$ is multiplicity-free by
Lemma~\ref{multfree}, it then follows that
$(Y_f+R)/R=Y(\sss,\eps)/R$. By the inductive hypothesis, we have
$R=\sum_{(\sss',\eps')\in\mathcal Z}\rad Y(\sss',\eps')\leq\rad
Y_<(\sss,\eps)$ and since $Y_<(\sss,\eps)/R$ is semisimple, we have
$R=\rad Y_<(\sss,\eps)$. We have therefore proved that $Y_f$ maps
onto $Y(\sss,\eps)/\rad Y_<(\sss,\eps)$. Then $Y_f$ contains a
submodule of $Y_<(\sss,\eps)$ which maps onto $Y_<(\sss,\eps)/\rad
Y_<(\sss,\eps)$. Since the radical of a module is the intersection
of the maximal submodules, the above submodule of $Y_<(\sss,\eps)$
must be all of $Y_<(\sss,\eps)$. Hence $Y_f$ contains
$Y_<(\sss,\eps)$ and we conclude that $Y_f=Y(\sss,\eps)$. The
theorem is proved.
\end{proof}

The following corollary is immediate.
\begin{corollary}\label{cyclicmodule}
Let $(\sss,\eps)\in \SSS[d]$. Then any $f\in
Y(\sss,\eps)[d]\setminus Y_<(\sss,\eps)[d]$ generates
$Y(\sss,\eps)[d]$.
\end{corollary}

The $k\GL(V)$-radical series of $Y(\sss)[d]$ ($[d]\neq[0]$)
and $Y(\sss)$ are given by digit sums as follows.
Let $\rad^i_{\GL(V)}M$ denote the $i$-th $k\GL(V)$-radical of the
$k\GL(V)$-module $M$. Then
\begin{equation}\label{GLrad}
\rad^i_{\GL(V)}Y(\sss)[d]=\sum_{\abs{\sss'}=\abs{\sss}-i}Y(\sss')[d],
\end{equation}
with a similar equation for $Y(\sss)$.
These results can be read off from \cite{BS}.

Our next result gives the analogous statements for $Y(\sss,\eps)[d]$
and $Y(\sss,\eps)$.
\begin{corollary}
\begin{enumerate}
\item[(i)] If $[d]\neq [0]$ then
$$
\rad^iY(\sss,\eps)[d]=\sum_{(\sss'',\eps'')}Y(\sss'',\eps'')[d]
$$
where the sum is over all $(\sss'',\eps'')\in\SSS[d]$ such that
$(\sss'',\eps'')\leq(\sss,\eps)$ and $\abs{\sss''}=\abs{\sss}-i$.
\item[(ii)] If $(\sss,\eps)\in\SSS$, then
$$
\rad^iY(\sss,\eps)=\sum_{(\sss'',\eps'')}Y(\sss'',\eps'')
$$
where the sum is over all $(\sss'',\eps'')\leq(\sss,\eps)$
such that $\abs{\sss''}=\abs{\sss}-i$.
\end{enumerate}
\end{corollary}
\begin{proof}
We will only prove (i), since (ii) is similar. Let $M_i$ denote the
module on the right side of the equation in (i). By
Remark~\ref{digitsum}, and Theorem~\ref{submodule} we see that
$M_{i+1}$ is the sum of all of the radicals of the
$Y(\sss',\eps')[d]$ occurring in $M_i$, and therefore $M_{i+1}\leq
\rad M_i$, since the radical of a sum of submodules of a module
contains the sum of their radicals. It remains to show that
$M_i/M_{i+1}$ is semisimple, which will show
$M_{i+1}\geq\rad M_i$, completing the proof.

We claim that
\begin{equation}
M_i=\rad^i_{\GL(V)}Y(\sss)[d]\cap Y(\sss,\eps)[d].
\end{equation}
From the claim, $M_i/M_{i+1}$ is isomorphic to a 
$k\Symp(V)$-submodule of the semisimple $k\GL(V)$-module
$(\rad^i_{\GL(V)}Y(\sss)[d])/(\rad^{i+1}_{\GL(V)}Y(\sss)[d])$;
so it is a semisimple $k\Symp(V)$-module, since every
simple $k\GL(V)$-composition factor is semisimple
as a $k\Symp(V)$-module.
To prove our claim, we consider the basis of $Y(\sss)[d]$
consisting of all symplectic basis functions with
$\HH$-types $\leq\sss$. The subset of this basis consisting
of those functions whose $\HH$-types
satisfy $\abs{\sss'}\leq\abs{\sss}-i$ form a basis of
$\rad^i_{\GL(V)}Y(\sss)[d]$, by the description of
$k\GL(V)$-radical series above.
By Lemma~\ref{signedbasis}, the subset of this basis consisting
of those elements whose signed types are $\leq (\sss,\eps)$
form a basis $Y(\sss,\eps)[d]$. And, by Lemma~\ref{signedbasis}
and the definition of $M_i$, the subset of the above basis
of $Y(\sss)[d]$ of functions whose signed $\HH$-types
satisfy both conditions $(\sss',\eps)\leq(\sss',\eps)$
and $\abs{\sss'}\leq\abs{\sss}-i$ form a basis of $M_i$.
The claim is established and the proof complete.
\end{proof}

\begin{remark}
The above corollary may be restated as saying the
$k\Symp(V)$-radical series of $Y(\sss,\eps)[d]$
and $Y(\sss,\eps)$ are given by intersecting the modules
with the $k\GL(V)$-radical series of $Y(\sss)[d]$ and
$Y(\sss)$, respectively.
\end{remark}

\section{The dimensions of ${\rm Im}(\eta_r)$}\label{dimensions}
Recall from Section~\ref{intro} that for $1\leq r\leq m$, $\eta_r$ denotes the
incidence map from $k[\Ii_r]$ to $k[P]$ sending a totally isotropic
$r$-dimensional subspace of $V$ to its characteristic function in
$P$.
For $m+1\leq r\leq 2m-1$, we can define $\Ii_r$ to be the set of
$r$-dimensional subspaces of the form $W^\perp=\{v\in V\mid \langle
v, w\rangle =0 , \text{ for all }w\in W\}$, for some totally
isotropic $(2m-r)$-dimensional subspace $W$. We can also consider
the incidence maps $\eta_r$ in this case.

\begin{theorem}\label{twoparts}
\begin{enumerate}
\item[(i)] We have
$${\rm Im}(\eta_m)=k1\oplus Y(\sss_m,\eps_m),$$
where $\sss_m=(m,m,\ldots ,m)$ and $\eps_m=\{0,1,\ldots ,t-1\}$.

\item[(ii)] If $1\leq r\leq 2m-1$ and $r\neq m$, then
$${\rm Im}(\eta_r)=k1\oplus Y(\sss_r),$$
where $\sss_r=(2m-r,2m-r,\ldots , 2m-r)$. In particular, if $1\leq
r< m$, then the $\Ff_q$-code generated by the characteristic
functions of all totally isotropic $r$-dimensional subspaces of $V$
is equal to the $\Ff_q$-code generated by the characteristic
functions of all $r$-dimensional subspaces of $V$.
\end{enumerate}
\end{theorem}

\begin{proof} We shall assume that $t>1$. When $t=1$, a similar and easier
argument works, but we omit the details to keep the argument clear,
since this case is already known \cite{sinsymp}.

(i) Since each point of $P$ is contained in
$\prod_{i=1}^{m-1}(1+q^i)$ totally isotropic $m$-dimensional
subspaces of $V$, by adding up the characteristic functions of all
totally isotropic $m$-dimensional subspaces of $V$, we get a nonzero
constant function. Hence $k1\subset {\rm Im}(\eta_m)$, where $k1$ is
the space of constant functions. Therefore we have a
$k\Symp(V)$-decomposition
$${\rm Im}(\eta_m)=k1\oplus M,$$
where $M\subset Y_P$ (cf. (\ref{splittingofk[P]})).

Let $L$ be the totally isotropic $m$-dimensional subspace of $V$
defined by the equations $x_i=0$, $i=1,2,\ldots ,m$, and $\chi_L$ be
the characteristic function of $L$. Since $\Symp(V)$ is transitive
on $\Ii_m$, we have
$${\rm Im}(\eta_m)=k\Symp(V)\chi_L.$$
Note that
\begin{eqnarray*}
\chi_L &=& (1-x_1^{q-1})(1-x_2^{q-1})\cdots (1-x_m^{q-1})\\
       &=& 1+f,
       \end{eqnarray*}
where $f=\sum_{\emptyset\neq I\subseteq \{1,2,\ldots
,m\}}(-1)^{|I|}\mathbf x_{I}^{q-1}$, and $\mathbf x_{I}$ stands for
$\prod_{i\in I}x_i$.  Therefore, we have $M=k\Symp(V)f$. For $0<|I|<m$,
the monomial $\mathbf x_{I}^{q-1}$ is a
symplectic basis function of signed type $((|I|,|I|,\ldots
,|I|),\emptyset)$, which lies below the signed
type $(\sss_m,\eps_m)$ of the symplectic basis function
$x_1^{q-1}x_2^{q-1}\cdots x_m^{q-1}$ in the poset $\SSS$. Hence $f\in
Y(\sss_m,\eps_m)\setminus Y_{<(\sss_m,\eps_m)}$. Therefore by
Corollary~\ref{cyclicmodule}, we have
$$M=Y(\sss_m,\eps_m).$$
We have proved (i).

(ii) First we deal with the case where $1\leq r <m$. Choose $L$ to be
the totally isotropic $r$-dimensional subspace of $V$ defined by the
equations $x_1=x_2=\cdots =x_m=0$ and $y_1=y_2=\cdots =y_{m-r}=0$.
Then the characteristic function of $L$ in $P$ is
$$\chi_L = (1-x_1^{q-1})(1-x_2^{q-1})\cdots
(1-x_m^{q-1})(1-y_1^{q-1})\cdots (1-y_{m-r}^{q-1}).$$ Since
$\Symp(V)$ is transitive on $\Ii_r$, we have ${\rm
Im}(\eta_r)=k\Symp(V)\chi_L$. This module also has the splitting
$$k\Symp(V)\chi_L=k 1\oplus N,$$
where $N=k\Symp(V)f$, $f=\chi_L -1$.  Note that
$$f=(-1)^rx_1^{q-1}\cdots x_m^{q-1}y_1^{q-1}\cdots y_{m-r}^{q-1} +
(-1)^{r-1}x_2^{q-1}\cdots x_m^{q-1}y_1^{q-1}\cdots
y_{m-r}^{q-1}+\cdots.$$ The symplectic basis function
$x_1^{q-1}\cdots x_m^{q-1}y_1^{q-1}\cdots y_{m-r}^{q-1}$ has signed
type $(\sss_r, \emptyset)$. The remaining terms in $f$ have signed
types strictly less than $(\sss_r, \emptyset)$. Hence by
Corollary~\ref{cyclicmodule}, we have $N=Y(\sss_r,\emptyset)$, which
in turn is equal to $Y(\sss_r)$ since $(\sss',\eps')\leq
(\sss_r,\emptyset)$ simply means $\sss'\leq \sss_r$. The proof of (ii) is complete in the case where $1\leq r<m$.
A similar argument works for the $m<r\leq 2m-1$ case.
\end{proof}

Next we develop the recursion for the $p$-ranks of the incidence
matrices between points and $m$-flats of $\Polar(2m-1, q)$ in terms of
$t$, where $q=p^t$, $p$ an odd prime. In particular, we will give a
proof for Theorem~\ref{prank}.

\begin{proposition}\label{summation}
Let $A_{1,m}^m(p^t)$ be the incidence matrix between points and
$m$-flats of $\Polar(2m-1, p^t)$, as defined in Section\ref{intro}. Assume that
$p$ is odd. Then
$${\rm rank}_p(A_{1,m}^m(p^t))= 1 + \sum_{\forall j, 1\leq s_j\leq
m}\prod_{j=0}^{t-1}d_{(s_j,s_{j+1})},$$ where
$$d_{(s_j,s_{j+1})} = \left\{\begin{array}{ll}
    {\rm dim}(S^+)=(d_{m(p-1)}+p^m)/2, & {\rm if}\ s_j=s_{j+1}=m, \\
    d_{\lambda_j},\ \mathrm{where}\ \lambda_j=ps_{j+1}-s_j, & \mathrm{otherwise.}
            \end{array}\right. $$
\end{proposition}

\begin{proof} By (i) of Theorem~\ref{twoparts}, the $p$-rank of
$A_{1,m}^m(p^t)$ is 1 plus the dimension of $Y(\sss_m,\eps_m)$,
where $\sss_m=(m,m,\ldots ,m)$ and $\eps_m=\{0,1,\ldots ,t-1\}$. By
Theorem~\ref{submodule}, the $k\Symp(V)$ module $Y(\sss_m,\eps_m)$
is multiplicity-free, and has as composition factors all
$L(\sss',\eps')$, $(\sss',\eps')\leq (\sss_m,\eps_m)$. Adding up the
dimensions of these composition factors (recall (\ref{dimen}) and
(\ref{dim+-})), we obtain the summation formula for ${\rm
rank}_p(A_{1,m}^m(p^t))$.
\end{proof}

\begin{corollary}
The $p$-rank of $A^m_{1,m}(p^t)$, when $p$ is an odd prime, is given by
$${\rm rank}_p(A_{1,m}^m(p^t))=1+\Trace(D^t)=1 + \alpha_1^t + \cdots +  \alpha_m^t,$$
 where
$$D=\begin{pmatrix}
d_{(1,1)} & d_{(1,2)}& \cdots & d_{(1,m)}\\
d_{(2,1)} & d_{(2,2)}& \cdots & d_{(2,m)}\\
 \vdots  & \vdots   & \ddots & \vdots\\
d_{(m,1)} & d_{(m,2)} & \cdots & d_{(m,m)}
\end{pmatrix},$$
and $\alpha_1, \alpha_2,\ldots , \alpha_m$ are the
eigenvalues of $D$.
\end{corollary}
Note that some of the entries of $D$ may be zero. We are now ready
to give the proof of Theorem~\ref{prank}.

\vspace{0.1in} \noindent {\bf Proof of Theorem~\ref{prank}:} We are
dealing with $\Polar(3,p^t)$, \emph{i.e.,} the case where $m=2$.  To
simplify notation, we write $A_{1,2}^2(p^t)$ simply as $A^{(t)}$. In
this case, we have
\begin{eqnarray*}
d_{(1,1)} &=&{\rm dim} (S^{p-1})=\frac {p(p+1)(p+2)} {6},\\
d_{(1,2)}&=&{\rm dim} (S^{2p-1})=\frac {2p(p-1)(p+1)} {3},\\
d_{(2,1)}&=&{\rm dim} (S^{p-2})=\frac {p(p+1)(p-1)} {6},\\
d_{(2,2)}&=&{\rm dim} (S^{+})=\frac {p(p+1)(2p+1)} {6}.
\end{eqnarray*}
Therefore
$$D=\frac {p(p+1)} {6}\begin{pmatrix}
p+2 & 4(p-1)\\
p-1 & 2p+1 \\
\end{pmatrix}.$$
This matrix $D$ has two distinct eigenvalues
$$\alpha_1,\alpha_2=\frac{p(p+1)^2}{4}\pm\frac{p(p+1)(p-1)}{12}\sqrt{17}.$$
Therefore we have
\begin{equation*}
{\rm rank}_p(A^{(t)})=1 +  \alpha_1^t +  \alpha_2^t.
\end{equation*}
\qed \vspace{0.1in}

The case where $m=3$ can be similarly analyzed. The matrix $D$ in
this case is given as follows.  
\begin{eqnarray*}
D&=&\frac{1}{120}
\begin{pmatrix}
(p+4)!/(p-1)! & (p^3-p)(p+2)(26p+48)
& 66p^5-210p^3+144p\\
(p+3)!/(p-2)!
 & 26p^5+50p^4+10p^3+10p^2+24p
& 66p^5-30p^3-36p  \\(p+2)!/(p-3)!
 &26p^5-10p^3-16p  & 33p^5+75p^3+12p \\
\end{pmatrix}
\end{eqnarray*}
The eigenvalues of $D$ have very complicated expressions: we will
not write down them here.

\vspace{0.2in}
\noindent{\bf Acknowledgements.} Machine computations for the case $q=9$ 
and the case $q=27$ done by Eric Moorhouse and Dave Saunders respectively 
were helpful in the early stages of our investigations.

\end{document}